\documentclass[reqno, a4paper]{amsart}

\usepackage[british]{babel}
\usepackage{amsmath}
\usepackage{amsfonts}
\usepackage{amssymb}
\usepackage{amsthm}
\usepackage[all]{xy}

\theoremstyle{plain}
\newtheorem{theorem}[subsection]{Theorem}
\newtheorem{lemma}[subsection]{Lemma}
\newtheorem{proposition}[subsection]{Proposition}

\theoremstyle{definition}
\newtheorem{definition}[subsection]{Definition}
\newtheorem{remark}[subsection]{Remark}
\newtheorem{notation}[subsection]{Notation}

\newcommand{\comp}{\circ}
\newcommand{\defn}{\textbf}
\newcommand{\To}{\Rightarrow}
\newcommand{\noproof}{\hfill \qed}
\newcommand{\meet}{\wedge}
\newcommand{\join}{\vee}

\newcommand{\coker}{\ensuremath{\mathrm{coker\,}}}
\renewcommand{\ker}{\ensuremath{\mathrm{ker\,}}}

\newcommand{\A}{\ensuremath{\mathcal{A}}}
\newcommand{\B}{\ensuremath{\mathcal{B}}}

\newcommand{\Smith}{\ensuremath{\mathrm{S}}}
\newcommand{\Huq}{\ensuremath{\mathrm{H}}}
\newcommand{\Higgins}{\ensuremath{\mathrm{\Omega}}}

\newcommand{\Ab}{\ensuremath{\mathsf{Ab}}}
\newcommand{\CExt}{\ensuremath{\mathsf{CExt}}}
\newcommand{\Ext}{\ensuremath{\mathsf{Ext}}}
\newcommand{\Arr}{\ensuremath{\mathsf{Arr}}}

\DeclareMathAlphabet{\mathitbf}{OML}{cmm}{b}{it}
\renewcommand{\a}{\mathitbf{a}}
\renewcommand{\b}{\mathitbf{b}}
\renewcommand{\k}{\mathitbf{k}}
\newcommand{\x}{\mathitbf{x}}
\newcommand{\m}{\mathitbf{m}}
\newcommand{\n}{\mathitbf{n}}
\newcommand{\p}{\mathitbf{p}}

\newbox\pullbackbox
\setbox\pullbackbox=\hbox{\xy 0;<1mm,0mm>: \POS(4,0)\ar@{-} (0,0) \ar@{-} (4,4)
\endxy}
\def\pullback{\copy\pullbackbox}

\newbox\pushoutbox
\setbox\pushoutbox=\hbox{\xy 0;<1mm,0mm>: \POS(0,4)\ar@{-} (0,0) \ar@{-} (4,4)
\endxy}
\def\pushout{\copy\pushoutbox}

\newdir{>}{{}*:(1,-.2)@^{>}*:(1,+.2)@_{>}}
\newdir{<}{{}*:(1,+.2)@^{<}*:(1,-.2)@_{<}}

\begin{document}

\title[Relative Commutator Theory in Semi-Abelian Categories]{Relative Commutator Theory\\ in Semi-Abelian Categories}

\author{Tomas Everaert}
\email{teveraer@vub.ac.be}
\address{Vakgroep Wiskunde, Vrije Universiteit Brussel, Pleinlaan 2, 1050 Brussel, Belgium}

\author{Tim Van~der Linden}
\email{tim.vanderlinden@uclouvain.be}
\address{Institut de recherche en math\'ematique et physique, Universit\'e catholique de Louvain, chemin du cyclotron~2, 1348 Louvain-la-Neuve, Belgium}
\address{Centro de Matem\'atica, Universidade de Coimbra, 3001--454 Coimbra, Portugal}

\thanks{The first author was supported by Fonds voor Wetenschappelijk Onderzoek (FWO-Vlaanderen). The second author works as \emph{charg\'e de recherches} for Fonds de la Recherche Scientifique--FNRS; his research was supported by Centro de Matem\'atica da Universidade de Coimbra and by Funda\c c\~ao para a Ci\^encia e a Tecnologia (grant number SFRH/BPD/38797/2007)}



\begin{abstract}
Basing ourselves on the concept of double central extension from categorical Galois theory, we study a notion of commutator which is defined relative to a Birkhoff subcategory $\B$ of a semi-abelian category~$\A$. This commutator characterises Janelidze and Kelly's $\B$-central extensions; when the subcategory $\B$ is determined by the abelian objects in $\A$, it coincides with Huq's commutator; and when the category $\A$ is a variety of $\Omega$-groups, it coincides with the relative commutator introduced by the first author. 
\end{abstract}

\maketitle

\section{Introduction}
The aim of this article is to fill in the question mark in the diagram 
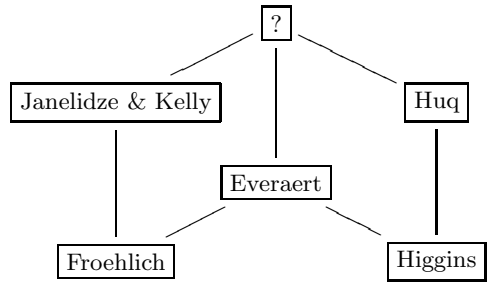
\[
\xymatrix@!0@C=6em@R=3em{& \fbox{\small ?} \ar@{-}[dl] \ar@{-}[dd] \ar@{-}[rd] \\
\fbox{\small Janelidze \& Kelly} \ar@{-}[dd] && \fbox{\txt{\small Huq}} \ar@{-}[dd]\\
& \fbox{\small Everaert} \ar@{-}[dl] \ar@{-}[rd] \\
\fbox{\txt{\small Froehlich}} && \fbox{\small Higgins}}
\]
which relates several non-equivalent concepts of \emph{commuting normal subobjects}, here named after the authors who introduced them. This diagram is meant to be read in the following manner. 

The bottom triangle restricts itself to theories which make sense for varieties of $\Omega$-groups, while the top triangle extends those theories to a categorical context. In the left hand side column we have theories which are \emph{one-dimensional} and \emph{relative}; the theories in the right hand side column, however, are \emph{two-dimensional} and \emph{absolute}, while the ones in the middle column are \emph{two-dimensional} and \emph{relative}. So we are looking for a \emph{categorical commutator theory which is both relative and two-dimensional}. Let us explain in more detail what this means for us.

\subsection{The bottom triangle}
Recall that a \defn{variety of $\Omega$-groups}~\cite{Higgins} is a variety in the sense of universal algebra which is pointed (i.e., it has exactly one constant) and has amongst its operations and identities those of the variety of groups. Apart from groups, the examples include the varieties of abelian groups, of non-unital rings, of commutative algebras, of modules and of Lie algebras, and also the categories of crossed modules and of precrossed modules are known (essentially from~\cite{LR}) to be equivalent to varieties of $\Omega$-groups.

In this context there are two classical approaches to commutator theory. On the one hand, there is the Higgins commutator of normal subobjects~\cite{Higgins} which has as particular cases the ordinary commutators of groups, rings, etc. It is \defn{two-dimensional} in the sense that any two normal subobjects (i.e., ideals or kernels) $N$ and $M$ of an object~$A$ in a variety of $\Omega$-groups~$\A$ have a commutator $[N,M]^{\Higgins}$, namely, the normal subobject of the join $M\join N=M\cdot N$ of $M$ and~$N$ generated by the set
\[
\{w(\m\n)w(\n)^{-1}w(\m)^{-1}\mid \text{$w$ is a term, $\m\in M$ and $\n\in N$}\},
\]
where the notation ``$\m\in M$'' means that $\m$ is a finite sequence $(m_{1},\dots,m_{r})$ of elements in $M$. Call an object~$A$ of $\A$ \defn{abelian} when it can be endowed with the structure of an internal abelian group (necessarily in a unique way). The subcategory of~$\A$ determined by the abelian objects is denoted by $\Ab\A$. It is well known (and easily verified) that when $\A$ is a variety of $\Omega$-groups, an algebra~$A$ is in~$\Ab\A$ precisely when the product map $A\times A\to A$ (sending a pair of elements $(a,a')$ to its product $aa'$) is a homomorphism in the variety. From this it follows immediately that the Higgins commutator characterises the abelian objects: $A$ is abelian if and only if $[A,A]^{\Higgins}=0$.

On the other hand there is the \defn{relative} notion of central extension due to Fr\"ohlich~\cite{Froehlich} (see also Lue~\cite{Lue} and Furtado-Coelho~\cite{Furtado-Coelho}). This notion of central extension corresponds to a \defn{one-dimensional} commutator. Here one starts from a variety of $\Omega$-groups $\A$ together with a chosen subvariety $\B$ of~$\A$. 

The subvariety $\B$ is completely determined by a set of identities of terms of the form $w(\x) = 1$; the set of all corresponding terms $w(\x)$ is denoted by
\[
W_{\B}=\{w(\x)\mid w(\b)=1,\forall B\in \B,\forall \b\in B\},
\]
and an object $A$ of $\A$ belongs to $\B$ if and only if $w(\a)=1$ for all $w\in W_{\B}$ and all~$\a\in A$.

An \defn{extension} $f\colon{A\to B}$ in $\A$ is a regular epimorphism, i.e., a surjective homomorphism. Let $K$ denote the kernel of $f$. The normal subobject $[K,A]^{\Omega}_{\B}$ of $A$ generated by the set
\[
\{w(\k\a)w(\a)^{-1}\mid \text{$w\in W_{\B}$, $\k\in K$ and $\a\in A$}\}
\]
is called the \defn{relative commutator (with respect to $\B$)} of $K$ and~$A$. (Note that Fr\"ohlich uses the notation $V_1$ for the relative commutator.) The extension~$f$ is \defn{central (with respect to $\B$)} when $[K,A]^{\Omega}_{\B}$ is zero. It is easily seen that this relative commutator characterises objects of~$\B$ as follows: $A$~belongs to~$\B$ if and only if $[A,A]^{\Omega}_{\B}$ is zero. 

In the \defn{absolute} case when the subvariety $\B$ consists of all abelian objects in $\A$, it was shown in~\cite{Furtado-Coelho} that the two commutators coincide,
\[
[K,A]^{\Omega}_{\Ab\A}=[K,A]^{\Higgins}.
\]
(Note here that $K\join A=A$.) The main advantage of the relative approach is that one may consider many situations which are not covered by the Higgins commutator. For instance, the notion of central extension of precrossed modules relative to the subvariety of crossed modules is of this type. The main advantage of the Higgins commutator is that it is two-dimensional. So the Higgins commutator is two-dimensional and absolute, the Fr\"ohlich commutator is one-dimensional and relative, and in the one-dimensional absolute case the two commutators coincide. What about the two-dimensional relative case? 

In his article~\cite{EveraertCommutator} the first author of the present article aims at answering precisely this question. He introduces a two-dimensional relative commutator for varieties of $\Omega$-groups which restricts to the Higgins commutator in the absolute case and which characterises Fr\"ohlich's relative central extensions. Given any pair of normal subobjects $M$ and $N$ of an object $A$ of $\A$, the commutator~$[M,N]_{\B}$ is the normal subobject of $M\join N$ generated by the set
\[
\{w(\m\n)w(\n)^{-1}w(\m)^{-1} w(\p)\mid w\in W_{\B}, \m\in M, \n\in N, \p\in M\meet N\}.
\]
The examples give an indication of how good his definition is. For instance, when considering the variety of precrossed modules together with the subvariety of crossed modules, the relative commutator obtained is the so-called \emph{Peiffer commutator}, which is exactly what one would expect.

\subsection{The left hand side column}
Basing themselves on ideas from categorical Galois theory~\cite{Janelidze:Pure, Borceux-Janelidze}, in the article~\cite{Janelidze-Kelly} Janelidze and Kelly introduce a general notion of central extension, relative with respect to a Birkhoff subcategory~$\B$ of a (Barr) exact category~$\A$. This notion of relative central extension is a generalisation of Fr\"ohlich's definition. 

In what follows, we shall restrict ourselves to the case of \defn{semi-abelian} categories~\cite{Janelidze-Marki-Tholen}: pointed, exact and protomodular with binary sums. So let~$\A$ be a semi-abelian category and $\B$ a \defn{Birkhoff subcategory} of $\A$---full, reflective and closed under subobjects and regular quotients; a Birkhoff subcategory of a variety is nothing but a subvariety. Let $I\colon{\A\to \B}$ denote the reflector, and $\eta\colon{1_{\A}\To I}$ the unit of the adjunction. Recall from~\cite{Janelidze-Kelly} that the closure of $\B$ under subobjects and regular quotients is equivalent to the condition that the commutative square 
\begin{equation}\label{Diagram-unitsquare}
\vcenter{\xymatrix@!0@=3.5em{A \ar[r]^-{f} \ar[d]_-{\eta_{A}} & B \ar[d]^-{\eta_{B}}\\
IA \ar[r]_-{If} & IB}}
\end{equation}
is a pushout of regular epimorphisms, for any regular epi $f\colon A\to B$.

An \defn{extension} in $\A$ is a regular epimorphism. Such an extension $f\colon{A\to B}$ is called \defn{trivial (with respect to $\B$)} when the induced commutative square~\eqref{Diagram-unitsquare} is a pullback. $f$ is \defn{central (with respect to $\B$)} when it is \emph{locally trivial} in the sense that there exists a regular epimorphism $p\colon E\to B$ such that the pullback $p^*(f)\colon E\times_B A\to E$ of $f$ along $p$ is a trivial extension. Since, in the present context, this implies that $f^*(f)$ is trivial, we have that $f$ is central if and only if it is \defn{normal}: either one of the projections in the kernel pair $(R[f],f_{0},f_{1})$ of~$f$ is a trivial extension. It is explained in the articles~\cite{Janelidze-Kelly, Bourn-Gran} why these central extensions reduce to Fr\"ohlich's when the category $\A$ is a variety of $\Omega$-groups.

This notion of relative central extension induces a \defn{one-dimensional relative commutator} as follows~\cite{EverVdL1,EGVdL}. Let $[-]_{\B}\colon{\A\to \A}$ denote the \emph{radical} induced by~$\B$: the functor which maps an object $A$ of $\A$ to the object $[A]_{\B}$ defined through the short exact sequence
\[
\xymatrix@!0@=3.5em{0 \ar[r] & [A]_{\B} \ar[r]^-{\mu_{A}} & A \ar[r]^-{\eta_{A}} & IA \ar[r] & 0,}
\]
and a morphism $a\colon{A'\to A}$ to its (co)restriction $[a]_{\B}\colon{[A']_{\B}\to [A]_{\B}}$. Let again $f\colon{A\to B}$ be an extension and let $K$ be its kernel. By protomodularity, $f$ is $\B$-central if and only if for the kernel pair $(R[f],f_{0},f_{1})$ of $f$, the (co)restrictions
\[
[f_{0}]_{\B},[f_{1}]_{\B}\colon{[R[f]]_{\B}\to [A]_{\B}}
\]
of the two projections are isomorphisms (see~\cite{Bourn-Gran}). Hence the kernel $[K,A]_{\B}$ of $[f_{0}]_{\B}$ measures how far $f$ is from being central: $f$ is $\B$-central if and only if~$[K,A]_{\B}$ is zero.

The object $[K,A]_{\B}$ may be considered as a normal subobject of $A$ via the composite
\[
\mu_{A}\comp [f_{1}]_{\B}\comp \ker [f_{0}]_{\B}\colon{[K,A]_{\B}\to A};
\]
the induced extension $A/[K,A]_{\B}\to B$ is the \defn{$\B$-centralisation} of $f$. We interpret~$[K,A]_{\B}$ as a commutator of $K$ with $A$, relative to the Birkhoff subcategory~$\B$ of $\A$. When $\A$ is a variety of $\Omega$-groups, $[K,A]_{\B}$ coincides with the relative commutator $[K,A]^{\Omega}_{\B}$, because they induce the same central extensions. And as in the varietal case, an object $A$ of $\A$ belongs to~$\B$ if and only if $[A,A]_{\B}=0$, because the extension $A\to 0$ is a split epimorphism, and therefore central if and only if it is trivial~\cite{Janelidze-Kelly}.

\subsection{The right hand side column}\label{Subsection-Huq-Commutator}
In his article~\cite{Huq}, Huq introduces a categorical notion of commutator of coterminal morphisms which makes sense in quite diverse algebraic settings. Using ``old-style'' axioms, he formulates his results for those categories we would nowadays call semi-abelian~\cite{Janelidze-Marki-Tholen}. Recast in more modern terminology by Bourn, his definition takes the following shape~\cite{Bourn-Huq}. In a semi-abelian category, consider two coterminal morphisms, $m\colon {M\to A}$ and $n\colon {N\to A}$, and the resulting square of solid arrows
\[
\xymatrix@!0@=3.5em{& M \ar[ld]_-{\langle 1_{M},0\rangle } \ar@{.>}[d] \ar[rd]^-{m}\\
M\times N \ar@{.>}[r] & Q & A. \ar@{.>}[l]|-{q}\\
& N \ar[lu]^-{\langle 0,1_{N}\rangle } \ar[ru]_-{n} \ar@{.>}[u]}
\]
The colimit of this square consists of an object $Q$ together with four morphisms with codomain $Q$ as indicated in the diagram. The morphism $q$ turns out to be a normal epimorphism; its kernel is denoted
\[
[m,n]^{\Huq}\colon {[M,N]^{\Huq}\to A}
\]
and called the \defn{Huq commutator of $m$ and~$n$}. It is convenient for us to restrict its use to the situation when $M$ and~$N$ are normal subobjects of $A$, i.e., $m$ and~$n$ are kernels. The commutator $[M,N]^{\Huq}$ becomes the ordinary commutator of normal subgroups $M$ and~$N$ in the case of groups, the ideal generated by $MN+NM$ in the case of rings, the Lie bracket in the case of Lie algebras, and so on. More generally, when computed in the join $M\join N$, we know from~\cite{Huq} that in any variety of $\Omega$-groups the Huq commutator~$[M,N]^{\Huq}$ coincides with the Higgins commutator $[M,N]^{\Higgins}$. Just as the Higgins commutator, the Huq commutator characterises the Birkhoff subcategory~$\Ab\A$ of $\A$ of abelian objects in $\A$. This is a consequence of the fact that, in a semi-abelian category $\A$, an object~$A$ admits at most one internal abelian group structure, and such a structure is entirely determined by a morphism $m\colon{A\times A\to A}$ which satisfies $m\comp \langle 1_{A},0\rangle = 1_{A}=m\comp \langle 0,1_{A}\rangle $~\cite{Huq,Bourn2002}.

\subsection{The question mark}\label{Question}
By now it is clear, we hope, that the purpose of the present article is to introduce a categorical version of the relative commutator for varieties of $\Omega$-groups, in such a way that
\begin{enumerate}
\item it characterises the $\B$-central extensions of $\A$,
\item it coincides with the Huq commutator when $\B$ is $\Ab\A$.
\end{enumerate}
In~\cite{EverVdL4} the present authors already introduced a relative concept of commuting normal subobjects, based on categorical Galois theory and valid in the context of semi-abelian categories. This notion was shown to be compatible with the relative commutator for varieties of $\Omega$-groups. What we still have to do now is
\begin{itemize}
\item explain how this induces a two-dimensional commutator;
\item prove that this commutator satisfies (1) and (2) above;
\item explore the commutator's basic properties.
\end{itemize}
One may ask whether it is worth the effort at all to leave the context of $\Omega$-groups and study a relative commutator from a categorical perspective. We claim that the categorical approach not only provides us with a conceptual explanation of the definitions (in terms of Galois theory) but also with interesting new examples. For instance, in the case of loops vs.\ groups considered in~\cite{EverVdL4}, the commutator becomes an \emph{associator}, and it effectively measures how well two normal subloops of a loop \emph{associate with each other}.

\subsection{Definition of the commutator} 
Let us now briefly sketch how the relative commutator $[-,-]_{\B}$ is defined. Let $\A$ again be a semi-abelian category and $\B$ a Birkhoff subcategory of $\A$. $M$ and~$N$ will be normal subobjects of an object $A$ of~$\A$. $R_{M}$ and~$R_{N}$ are the equivalence relations on the join $M\join N$ (taken in the lattice of normal subobjects of $A$) corresponding to $M$ and~$N$, and 
\[
\xymatrix@!0@R=3.5em@C=6em{R_{M}\square R_{N} \ar@<.5ex>[r]^-{r_{1}} \ar@<-.5ex>[r]_-{r_{0}} \ar@<.5ex>[d]^-{p_{1}} \ar@<-.5ex>[d]_-{p_{0}} & R_{N} \ar@<.5ex>[d] \ar@<-.5ex>[d]\\
R_{M} \ar@<.5ex>[r] \ar@<-.5ex>[r] & M\join N}
\] 
is the largest double equivalence relation on $R_{M}$ and $R_{N}$: the object $R_{M}\square R_{N}$ ``consists of'' all quadruples $(x, y, z, t)\in M\join N$ where $(x,z)$, $(y,t)\in R_{M}$ and $(x,y)$, $(z,t)\in R_{N}$.

The commutator of $M$ and $N$ is the meet
\[
[M,N]_{\B}=K[[p_{0}]_{\B}]\meet K[[r_{0}]_{\B}]
\]
of the kernels of the morphisms $[p_{0}]_{\B}$ and $[r_{0}]_{\B}$ in the following diagram, obtained by applying the functor $[-]_{\B}$ to the diagram above.
\begin{equation}\label{Commutator-Squares}
\vcenter{\xymatrix@!0@R=3.5em@C=6em{[R_{M}\square R_{N}]_{\B} \ar@<.5ex>[r]^-{[r_{1}]_{\B}} \ar@<-.5ex>[r]_-{[r_{0}]_{\B}} \ar@<.5ex>[d]^-{[p_{1}]_{\B}} \ar@<-.5ex>[d]_-{[p_{0}]_{\B}} & [R_{N}]_{\B} \ar@<.5ex>[d] \ar@<-.5ex>[d]\\
[R_{M}]_{\B} \ar@<.5ex>[r] \ar@<-.5ex>[r] & [M\join N]_{\B}}}
\end{equation}
It may be considered as a normal subobject of $M\join N$.

\subsection{Interpretation in terms of double central extensions}\label{Galois-Structure} 
We have to explain why $[M,N]_{\B}$ is defined the way it is. The reason comes from categorical Galois theory, in particular the theory of \emph{higher central extensions}. Just like the concept of central extension which is defined with respect to the adjunction
\begin{equation}\label{Adjunction-I}
\xymatrix@!0@=3.5em{{\A} \ar@<1ex>[r]^-{I} \ar@{}[r]|-{\perp} & {\B,} \ar@<1ex>[l]^-{\supset}}
\end{equation}
one may consider double central extensions which are defined with respect to the reflection of extensions to central extensions---the adjunction
\begin{equation}\label{Adjunction-I_{1}}
\xymatrix@!0@=6em{{\Ext\A} \ar@<1ex>[r]^-{I_{1}} \ar@{}[r]|-{\perp} & {\CExt_{\B}\A} \ar@<1ex>[l]^-{\supset}}
\end{equation}
where $\Ext\A$ is the category of extensions and commutative squares between them, and $\CExt_{\B}\A$ its full subcategory determined by those extensions which are central. The reflector $I_{1}$ takes an extension $f\colon{A\to B}$ with kernel $K$ and maps it to the central extension
\[
I_{1}f\colon{A/[K,A]_{\B}\to B}.
\]
This may be repeated ad infinitum, so that notions of \emph{$n$-fold central extension} are obtained, but for the present purposes the second step is sufficient. Double central extensions, first introduced by Janelidze for groups~\cite{Janelidze:Double}, are an important tool in semi-abelian (co)homology~\cite{EGVdL, Janelidze:Hopf-talk, RVdL}, and turn out to be precisely what is needed to understand how the relative commutator works. We refer the reader to the articles~\cite{EGVdL,EverHopf} for more details on higher central extensions. 

As we explain below, the commutator $[M,N]_{\B}$ is zero if and only if any (hence, all) of the four commutative squares in the diagram~\eqref{Commutator-Squares} is a pullback. Galois theory shows that this condition is equivalent to the square
\begin{equation}\label{Double-Central-Extension}
\vcenter{\xymatrix@!0@R=3.5em@C=6em{M\join N \ar[r]^-{q_{M}} \ar[d]_-{q_{N}} & \tfrac{M\join N}{M} \ar[d]\\
\tfrac{M\join N}{N} \ar[r] & 0}}
\end{equation}
being a double central extension. (Here $q_M$ denotes the cokernel of the normal monomorphism $M\to M\join N$.) When this happens, we say that $M$ and $N$ \defn{commute (with respect to $\B$)}.

Accordingly, given any two normal subobjects $M$ and $N$ of an object~$A$, the commutator $[M,N]_{\B}$ is the smallest normal subobject $J$ of $M\join N$ such that $M/J$ and $N/J$ commute; it is the normal subobject which must be divided out of $M\join N$ to turn the double extension~\eqref{Double-Central-Extension} into a double central extension.

\subsection{Structure of the text} In the following sections we shall explain why the commutator has the properties (1) and (2) mentioned in~\ref{Question}. With this purpose in mind, the text is structured as follows. In Section~\ref{Section-Preliminaries} we provide the necessary background for understanding the definition of the commutator: semi-abelian categories, normal subobjects, double extensions and double central extensions. Its basic technical properties and the proof of (1) are given in Section~\ref{Section-Definition}. In Section~\ref{Section-Huq} we prove~(2): the commutator $[-,-]_{\B}$ coincides with the Huq commutator in case~$\B$ is~$\Ab\A$. Finally, Section~\ref{Further-Remarks} brings up some further remarks and unanswered questions.

\tableofcontents
\pagebreak 

\section{Preliminaries}\label{Section-Preliminaries}
We recall some basic definitions and results which we shall need in the following sections.

\subsection{Semi-abelian categories}
A category is \defn{regular} when it is finitely complete with coequalisers of kernel pairs and with pullback-stable regular epimorphisms~\cite{Barr}. In a regular category, any morphism $f$ may be factored as a regular epimorphism followed by a monomorphism (called the \defn{image} of $f$), and this \defn{image factorisation} is unique up to isomorphism. Given a monomorphism $m\colon {M\to A}$ and a regular epimorphism $f\colon A\to B$, the \defn{direct image} $f(m)\colon{fM\to B}$ of $m$ along $f$ is the image of the composite~${f\comp m}$.

When a category is pointed and regular, \defn{protomodularity} can be defined via the following property, which is equivalent to the Short Five Lemma~\cite{Bourn1991,Bourn2001}: given any commutative diagram
\begin{equation}\label{Short-Five-Lemma}
\vcenter{\xymatrix@!0@=4em{K[f'] \ar[r]^-{\ker f'} \ar[d]_-k & A' \ar[r]^-{f'} \ar[d]_-a & B' \ar[d]^-b \\ K[f] \ar[r]_-{\ker f} & A \ar[r]_-{f} & B}}
\end{equation}
such that $f$ and $f'$ are regular epimorphisms, $k$ is an isomorphism if and only if the right hand square $b\comp f'=f\comp a$ is a pullback. (Here, we use the notation $\ker f\colon K[f]\to A$ for the kernel of $f$.) A \defn{homological} category is pointed, regular and protomodular~\cite{Borceux-Bourn}. In such a category, a regular epimorphism is always the cokernel of its kernel, and there is the following notion of short exact sequence. A \defn{short exact sequence} is any sequence
\[
\xymatrix@!0@=3.5em{K \ar[r]^-{k} & A \ar[r]^-{f} & B}
\]
with $k=\ker f$ and $f$ a regular epimorphism. We denote this situation by
\[
\xymatrix@!0@=3.5em{0 \ar[r] & K \ar[r]^-{k} & A \ar[r]^-{f} & B \ar[r] & 0.}
\] 
The following property holds.
\begin{lemma}\cite{Bourn2001}\label{Lemma-Pullback}
Consider a morphism of short exact sequences such as~\eqref{Short-Five-Lemma} above. The left hand side square $\ker f\comp k=a\comp \ker f'$ is a pullback if and only if $b$ is a mono.
\noproof
\end{lemma}

A \defn{(Barr) exact} category is regular and such that every internal equivalence relation is a kernel pair~\cite{Barr}. A homological category is exact if and only if the direct image of a normal monomorphism along a regular epimorphism is again a normal monomorphism. A \defn{semi-abelian} category is homological and exact with binary coproducts~\cite{Janelidze-Marki-Tholen}.

A \defn{regular pushout square} is a commutative square
\begin{equation}\label{Diagram-Square}
\vcenter{\xymatrix@!0@=3.5em{X \ar[r]^-{c} \ar[d]_-{d} & C \ar[d]^-{g}\\
D \ar[r]_-{f} & Z}}
\end{equation}
such that all its maps and the comparison map $\langle d,c\rangle \colon{X\to D\times_Z C}$ to the pullback of $f$ with $g$ are regular epimorphisms. In a semi-abelian category, every pushout of a regular epimorphism along a regular epimorphism is a regular pushout~\cite{Carboni-Kelly-Pedicchio}, and the following dual to Lemma~\ref{Lemma-Pullback} holds:
\begin{lemma}\label{Rotlemma}\cite{Bourn-Gran}
Given a morphism of short exact sequences such as~\eqref{Short-Five-Lemma} above with $a$ and $b$ regular epi, the right hand side square $f\comp a=b\comp f'$ is a (regular) pushout if and only if $k$ is a regular epimorphism. \noproof
\end{lemma}

\subsection{Normal subobjects}
A \defn{normal subobject} $N$ of an object~$A$ of a semi-abelian category is a subobject represented by a normal monomorphism $n\colon{N\to A}$. Let $M$ and~$N$ be two normal subobjects of $A$ with representing normal monomorphisms $m$ and~$n$. Taking into account Lemma~\ref{Lemma-Pullback} and the stability of normal monomorphisms under regular images, we may always form the $3\times 3$ diagram in Figure~\ref{3x3-Diagram} (in which all rows and columns are short exact sequences). The meet ${M\meet N}$ and the join $M\join N$ of the subobjects~$M$ and~$N$ are taken in the lattice of normal subobjects of $A$. We see that ${M\meet N}$ is computed as the pullback~\texttt{(i)} and~$M\join N$ is obtained through the pushout~\texttt{(ii)}, as the kernel of the composite morphism ${A\to {A}/{(M\join N)}}$. Of course, $M\meet N$ coincides with the meet~$M\cap N$ in the lattice of (all) subobjects of~$A$. One could also compute the join of~$M$ and~$N$ as (ordinary) subobjects of $A$ by taking the image $M\cup N$ of the morphism~$\langle \begin{smallmatrix}m \\ n\end{smallmatrix}\rangle\colon {M+N\to A}$. It is known~\cite{Borceux-Semiab,Huq} that both constructions yield the same result. We shall give an alternative proof of this fact below, but first we prove a weaker property. 

\begin{figure}
\[
\xymatrix@!0@=3.5em{& 0 \ar[d] & 0 \ar[d] & 0 \ar[d]\\
0 \ar[r] & M\meet N \ar@{}[rd]|-{\texttt{(i)}} \ar[r] \ar[d] \ar@{}[rd]|<<{\pullback} & N \ar[d]^-{n} \ar[r] & \tfrac{N}{M\meet N} \ar[d] \ar[r] & 0\\
0 \ar[r] &M \ar[r]_-{m} \ar[d] & A \ar[r] \ar[d] \ar@{}[rd]|>{\pushout} \ar@{}[rd]|-{\texttt{(ii)}} & \tfrac{A}{M} \ar[d] \ar[r] & 0\\
0 \ar[r] & \tfrac{M}{M\meet N} \ar[r] \ar[d] & \tfrac{A}{N} \ar[r] \ar[d] & \tfrac{A}{M\join N} \ar[d] \ar[r] & 0\\
& 0 & 0 & 0 }
\]
\caption{The $3\times 3$ diagram induced by $M$, $N$ normal in $A$}\label{3x3-Diagram}
\end{figure}

Let us fix some notation: we write $j$ for the normal monomorphism representing~$M\join N$, and $m'\colon M\to M\join N$ and $n'\colon N\to M\join N$ for the induced factorisations. Since $m'$ and $n'$ are normal monomorphisms, we may also consider the join of $M$ and $N$ as normal subobjects of $M\join N$. We denote it by $M\curlyvee N$ and write $j'\colon M\curlyvee N\to M\join N$ for the representing normal monomorphism. 
\begin{lemma}
The two joins $M\join N$ and $M\curlyvee N$ coincide: $j'$ is an isomorphism.
\end{lemma}
\begin{proof}
First of all note that the commutative square
\[
\xymatrix@!0@R=3.5em@C=5em{
{M\join N} \ar[r] \ar[d]_j & \frac{M\join N}{M} \ar[d]\\
A \ar[r] & \frac{A}{M}}
\]
is a pullback by protomodularity, so that the right hand vertical morphism is a monomorphism because, in a protomodular category, pullbacks reflect monos~\cite{Bourn1991}. (One could, alternatively, use Lemma~\ref{Lemma-Pullback} to prove that this morphism is a monomorphism.) Now, the normal monomorphisms~$m'$ and~$n'$ induce a $3\times 3$ diagram similar to Figure~\ref{3x3-Diagram}, and $j$ induces a morphism between the two $3\times 3$ diagrams, of which we consider only the last row: 
\[
\xymatrix@!0@R=3.5em@C=5em{
0 \ar[r] & \frac{N}{M\meet N} \ar[r] \ar@{=}[d] & \frac{M\join N}{M} \ar[r] \ar[d] & \frac{M\join N}{M\curlyvee N} \ar[d] \ar[r] & 0\\
0 \ar[r] & \frac{N}{M\meet N} \ar[r] & \frac{A}{M} \ar[r] & \frac{A}{M\join N} \ar[r] & 0}
\]
We have just explained why the middle vertical morphism is a monomorphism. Hence, using the same arguments as above, we find that also the right hand vertical morphisms is a mono. Since the composite 
\[
M\join N \to (M\join N)/(M\curlyvee N)\to A/(M\join N)
\]
is zero, we find that $(M\join N)/(M\curlyvee N)=0$, i.e., the factorisation $j'$ is an isomorphism.
\end{proof}

Now, taking this lemma into account, when $A=M\join N$ in the $3\times 3$ diagram above, the object $A/(M\join N)$ is zero, and we regain the \defn{Noether isomorphisms}~\cite{Borceux-Bourn}
\begin{equation}\label{Noether}
\frac{N}{M\meet N}\cong \frac{M\join N}{M}\qquad\text{and}\qquad\frac{M}{M\meet N}\cong \frac{M\join N}{N}.
\end{equation}
We are ready to prove the identity $M\join N=M\cup N$.

\begin{notation}\label{Notation-q}
Given a normal subobject $J$ of an object~$A$, the induced quotient of $A$ is denoted
\[
q_{J}^{A}\colon{A\to A/J};
\]
we write $R^{A}_{J}$ for the kernel pair $A\times_{A/J} A$ of $q^{A}_{J}$. 

Most of the time $A$ will be a join $M\vee N$, in which case we drop the $A$ from the notation and simply write 
\[
q_{J}\colon{M\join N\to (M\join N) / J}
\]
for the quotient and $R_{J}$ for the kernel pair of $q_{J}$. 
\end{notation}

\begin{proposition}\cite{Borceux-Semiab,Huq}\label{Lemma-Join}
If $M$ and $N$ are normal in $A$, then their join as normal subobjects $M\join N$ coincides with their join as subobjects $M\cup N$. Hence the morphism
\[
\langle \coker n,\coker m\rangle \colon{A\to (A/N)\times (A/M)}
\]
is a regular epimorphism if and only if such is the morphism
\[
\langle\begin{smallmatrix}m \\ n\end{smallmatrix}\rangle\colon {M+N\to A}.
\]
\end{proposition}
\begin{proof}
If $J$ is a subobject of $M\join N$ containing $M$ and $N$, then by Lemma~\ref{Lemma-Pullback} it induces a factorisation of the first of the isomorphisms~\eqref{Noether} as a morphism ${N/(M\meet N)\to J/M}$ followed by a monomorphism
\[
j\colon{J/M\to {(M\join N)}/{M}}.
\]
This $j$ is also a split epimorphism; hence it is an isomorphism, and $J$ is equal to~$M\join N$ by the Short Five Lemma. 

Now $M\cup N$ is a subobject of $M\join N$ containing $M$ and $N$, and the two joins coincide. 

As to the latter statement, the first condition holds if and only if the square
\[
\xymatrix@!0@R=3.5em@C=5em{A \ar[r] \ar[d] & \tfrac{A}{M} \ar[d]\\
\tfrac{A}{N} \ar[r] & 0}
\]
is a regular pushout. Since, in a semi-abelian category, a pushout of regular epimorphisms is necessarily regular, this happens when $A=M\join N$. But then~$A$ is~$M\cup N$ by the former part of the proof, and the second condition holds only when this is the case.
\end{proof}

Given a monomorphism $m\colon{M\to A}$, the \defn{normal closure} $\overline{M}^{A}$ of $M$ in $A$ always exists, and is computed as the kernel of the cokernel of $m$. It is the smallest normal subobject of $A$ that contains~$M$.

\subsection{Double (central) extensions}\label{Subsection-Double-Central-Extensions}
A \defn{double extension} is a regular pushout square~\eqref{Diagram-Square}. For instance, given any two normal subobjects $M$ and~$N$ of an object $A$ of $\A$, the induced pushout square~\eqref{Double-Central-Extension} is a double extension. Recall from~\cite{EGVdL} that pullbacks of double extensions exist in $\Ext\A$ and are degree-wise pullbacks in $\A$. Moreover, double extensions are pullback-stable. The category of double extensions in $\A$ and commutative cubes between them is denoted~$\Ext^{2}\!\A$. 

Double central extensions are defined with respect to the adjunction~\eqref{Adjunction-I_{1}} in the same way as central extensions are defined with respect to the adjunction~\eqref{Adjunction-I}. More precisely, a double extension~\eqref{Diagram-Square}, considered as a map $(c,f)\colon{d\to g}$ in the category $\Ext \A$, is \defn{trivial} when the left hand commutative square below, induced by the unit of the adjunction~\eqref{Adjunction-I_{1}}, is a pullback in~$\Ext\A$; this means that the right hand commutative square, in which the vertical morphisms are the canonical quotient maps, is a pullback in~$\A$.
\[
\vcenter{\xymatrix@1@!0@=3.5em{
d \ar[r]^{(c,f)} \ar[d] & g \ar[d]\\
I_1d \ar[r] & I_1g}}
\qquad\qquad
\vcenter{\xymatrix@!0@R=3.5em@C=5em{
X \ar[r]^c \ar[d] & C \ar[d] \\
\frac{X}{[K[d],X]_{\B}} \ar[r] & \frac{C}{[K[g],C]_{\B}}}}
\]
The square~\eqref{Diagram-Square} is a \defn{double central extension (with respect to $\B$)} when its pullback along some double extension is a trivial double extension. It is a \defn{double normal extension (with respect to $\B$)} when the first projection of its kernel pair
\[
\vcenter{\xymatrix@!0@R=3.5em@C=5em{R[c] \ar[r]^-{c_{0}} \ar[d]_-{r} & X \ar[d]^-{d} \\
R[f] \ar[r]_-{f_{0}} & D}}
\]
is a trivial double extension. (Alternatively, one could use the square of second projections.) By protomodularity, this amounts to the (one-dimen\-sional, relative) commutators $[K[r],R[c]]_{\B}$ and $[K[d],X]_{\B}$ being isomorphic. Similar to the one-dimensional case, double central extensions and double normal extensions coincide.

\subsection{Higher extensions} 
In what follows we shall also need three-fold extensions, so let us recall the definition of $n$-fold extension for arbitrary $n$. Given $n\geq 0$, denote by $\Arr^n\!\A$ the category of $n$-dimensional arrows in $\A$. (Zero-dimensional arrows---as well as zero-dimensional extensions---are just objects of $\A$.) A \defn{(one-fold) extension} is a regular epimorphism in $\A$. For $n\geq 1$, an \defn{$(n+1)$-fold extension} is a commutative square~\eqref{Diagram-Square} in~$\Arr^{n-1}\!\A$ (an arrow in~$\Arr^n\!\A$) such that all its maps and the comparison map $\langle d,c\rangle\colon{X\to D\times_Z C}$ to the pullback of~$f$ with~$g$ are~$n$-fold extensions. Thus for $n=2$ we regain the notion of double extension. 

A three-fold extension is a commutative cube
\[
\vcenter{\xymatrix@!0@=3.5em{& X \ar[rr] \ar@{.>}[dd] && C \ar[dd] \\
X' \ar[rr] \ar[dd] \ar[ru] && C' \ar[dd] \ar[ru]\\
& D \ar@{.>}[rr] && Z\\
D' \ar[rr] \ar@{.>}[ru] && Z' \ar[ru]}}
\qquad\qquad
\vcenter{\xymatrix@!0@R=3.5em@C=4.5em{X' \ar[r] \ar[d] & D'\times_{Z'}C' \ar[d] \\
X \ar[r] & D\times _{Z}C}}
\]
of which all faces as well as the induced right-hand square are double extensions. Since, in a semi-abelian category, regular epimorphisms are normal, the three-fold extension above is completely determined by the object $X'$ and the three normal subobjects given by the kernels of its ``initial ribs'' ${X'\to X}$, ${X'\to C'}$ and ${X'\to D'}$. Conversely, given an object $X'$ and three normal subobjects $J$,~$M$ and $N$ of $X'$, the following lemma determines when the induced cube is a three-fold extension.

\begin{lemma}\label{Lemma-Three-Fold-Extension}
Given normal subobjects $J$,~$M$ and $N$ of an object $X'$ in a semi-abelian category, the cube obtained by pushing out the induced quotients is a three-fold extension if and only if
\[
{q^{X'}_{J}(M\meet N)=q^{X'}_{J}M\meet q^{X'}_{J}N}.
\]
\end{lemma}
\begin{proof}
Since, in a semi-abelian category, pushouts of regular epimorphisms are regular, the induced cube is a three-fold extension as soon as the square 
\[
\vcenter{\xymatrix@!0@R=4em@C=7em{X' \ar[r] \ar[d] & \tfrac{X'}{M}\times_{\tfrac{X'}{M\join N}}\tfrac{X'}{N} \ar[d] \\
\tfrac{X'}{J} \ar[r] & \tfrac{X'}{J\join M}\times_{\tfrac{X'}{J\join M\join N}}\tfrac{X'}{J\join N}}}
\]
is a double extension. We already know that all morphisms in this square are regular epimorphisms, so by Lemma~\ref{Rotlemma} it is a double extension if and only if $q^{X'}_{J}(M\meet N)=q^{X'}_{J}M\meet q^{X'}_{J}N$.
\end{proof}

Further results on higher-dimensional extensions and central extensions may be found in~\cite{EGVdL} and~\cite{EverHopf}. Let us just recall here that, for any $n\geq 0$, a split epimorphism of $n$-fold extensions is always an $(n+1)$-fold extension, and it is an $(n+1)$-fold central extension if and only if it is a trivial $(n+1)$-fold extension.

Higher-dimensional central extensions are important in homology where they appear in the higher Hopf formulae, and in cohomology where (in the absolute case, and in low dimensions) they are classified by the cohomology groups~\cite{Gran-VdL,RVdL}.

\section{Definition and basic properties}\label{Section-Definition}
In this section we recall the categorical definition of the relative commutator from the introduction and we explore its basic properties: compatibility with the central extensions introduced by Janelidze and Kelly (Proposition~\ref{Proposition-Characterisation-Central-Extensions}), basic stability properties (Theorem~\ref{Proposition-Commutator-Properties}) and the case of $\Omega$-groups (Proposition~\ref{Theorem-Characterisation}).

In what follows, $\A$ will be a semi-abelian category and $\B$ a Birkhoff subcategory of $\A$.

\begin{definition}\label{Definition-Commuting-Subobjects}
Let $M$ and $N$ be normal subobjects of an object $A$ of $\A$. We say that $M$ and $N$ \defn{commute (with respect to $\B$)} when the double extension
\begin{equation}\label{Double-Extension-MN}
\vcenter{\xymatrix@!0@R=3.5em@C=5em{M\join N \ar[r]^-{q_{M}} \ar[d]_-{q_{N}} & \tfrac{M\join N}{M} \ar[d]\\
\tfrac{M\join N}{N} \ar[r] & 0}}
\end{equation}
is central (with respect to $\B$).
\end{definition}

Is is immediately clear that this notion of commuting subobjects characterises the $\B$-central extensions of $\A$ and the objects of $\B$: 

\begin{proposition}\label{Proposition-Characterisation-Central-Extensions}
An extension $f\colon{A\to B}$ in $\A$ is $\B$-central if and only if the object $A$ and the kernel $K$ of $f$ commute. An object $A$ of $\A$ lies in $\B$ if and only if~$A$ commutes with itself.
\end{proposition}
\begin{proof}
The first result holds because the double extension
\[
\xymatrix@!0@=3.5em{A \ar[r]^-{q_{A}} \ar[d]_-{f=q_{K}} & 0 \ar@{=}[d] \\
B \ar[r] & 0,}
\]
being a split epimorphism of extensions, is central if and only if it is trivial, which happens precisely when $f$ is a central extension. The second result follows from the first, since $A$ is in $\B$ if and only if the split epimorphism ${A\to 0}$ is a $\B$-central extension.
\end{proof}

\begin{lemma}\cite[Proposition 2.9]{EverVdL4}\label{Lemma-Characterisation-Square}
Let $M$ and $N$ be normal subobjects of an object~$A$. $M$ and $N$ commute if and only if any of the four commutative squares in the diagram 
\begin{equation}\label{Square-Commutator}
\vcenter{\xymatrix@!0@C=7em@R=3.5em{[R_{M}\square R_{N}]_{\B} \ar@<.5ex>[r]^-{[r_{1}]_{\B}} \ar@<-.5ex>[r]_-{[r_{0}]_{\B}} \ar@<.5ex>[d]^-{[p_{1}]_{\B}} \ar@<-.5ex>[d]_-{[p_{0}]_{\B}} & [R_{N}]_{\B} \ar@<.5ex>[d]^-{[\pi_{1}]_{\B}} \ar@<-.5ex>[d]_-{[\pi_{0}]_{\B}}\\
[R_{M}]_{\B} \ar@<.5ex>[r]^-{[\rho_{1}]_{\B}} \ar@<-.5ex>[r]_-{[\rho_{0}]_{\B}} & [M\join N]_{\B}}}
\end{equation}
is a pullback.\noproof
\end{lemma}

\begin{definition}\label{Definition-Commutator}
Let $M$ and $N$ be normal subobjects of an object $A$. Let
\[
[R_{M}]_{\B}\times_{[M\join N]_{\B}}[R_{N}]_{\B}
\]
denote the pullback of the morphisms $[\pi_{0}]_{\B}$ and $[\rho_{0}]_{\B}$ from Diagram~\eqref{Square-Commutator}. The \defn{commutator $[M,N]_{\B}$} is the kernel of the morphism
\[
\langle [p_{0}]_{\B},[r_{0}]_{\B}\rangle \colon{[R_{M}\square R_{N}]_{\B}\to [R_{M}]_{\B}\times_{[M\join N]_{\B}}[R_{N}]_{\B}},
\]
considered as a normal subobject of $M\join N$.
\end{definition}

\begin{remark}
Two normal subobjects $M$ and $N$ of an object $A$ commute if and only if $[M,N]_{\B}$ is zero. Indeed, the morphism $\langle [p_{0}]_{\B},[r_{0}]_{\B}\rangle $ is a regular (hence, normal) epimorphism because the square $[\pi_0]_{\B}\comp [r_0]_{\B}=[\rho_0]_{\B}\comp [p_0]_{\B}$ is a double extension as a split epimorphism of split epimorphisms. Hence its kernel is zero if and only if it is an isomorphism---which, by Lemma~\ref{Lemma-Characterisation-Square}, means that~$M$ and~$N$ commute. 
\end{remark}

\begin{remark}\label{Remark-Subobject}
The kernel of $\langle [p_{0}]_{\B},[r_{0}]_{\B}\rangle $ may indeed be considered as a normal subobject of $M\join N$, namely, through the composition of
\[
\ker \langle [p_{0}]_{\B},[r_{0}]_{\B}\rangle \colon{K[\langle [p_{0}]_{\B},[r_{0}]_{\B}\rangle ]\to [R_{M}\square R_{N}]_{\B}}
\]
with
\[
\rho_{1}\comp p_{1}\comp \mu_{R_{M}\square R_{N}}\colon{[R_{M}\square R_{N}]_{\B}\to M\join N.}
\]
First of all, this composite is a monomorphism, because
\[
\mu_{R_{M}\square R_{N}}\comp \ker \langle [p_{0}]_{\B},[r_{0}]_{\B}\rangle= \ker \langle p_0,r_0 \rangle \comp \mu_{K[\langle [p_{0}]_{\B},[r_{0}]_{\B}\rangle]}
\]
and both $\mu_{K[\langle [p_{0}]_{\B},[r_{0}]_{\B}\rangle]}$ and $\rho_{1}\comp p_{1}\comp \ker \langle p_0,r_0 \rangle$ are monomorphisms.

Now $\mu_{R_{M}\square R_{N}}\comp \ker \langle [p_{0}]_{\B},[r_{0}]_{\B}\rangle$ is a normal monomorphism as a meet of two normal monomorphisms. This follows from Lemma~\ref{Lemma-Pullback}, since the induced morphism
\[
{[R_{M}]_{\B}\times_{[M\join N]_{\B}}[R_{N}]_{\B}\to R_{M}\times_{M\join N}R_{N}}
\]
is a monomorphism. Hence 
\[
\rho_{1}\comp p_{1}\comp \mu_{R_{M}\square R_{N}}\comp \ker \langle [p_{0}]_{\B},[r_{0}]_{\B}\rangle
\]
is normal, being the direct image of this latter normal monomorphism along the regular epimorphism $\rho_{1}\comp p_{1}$.
\end{remark}

\begin{remark}
On the other hand, there is no reason why $[M,N]_{\B}$ should be a normal subobject of $A$. A counterexample is given in~\cite{MM-NC}.
\end{remark}

\begin{remark}
The commutator $[M,N]_{\B}$ is nothing but $L_{2}$ of the double extension~\eqref{Double-Extension-MN} as considered in the article~\cite{EGVdL}. 
\end{remark}

\begin{theorem}\label{Proposition-Commutator-Properties}
Let $M$, $N$, $L$ (resp.~$M'$, $N'$) be normal subobjects of an object~$A$ (resp.~$A'$). Let $J$ be a normal subobject of $M\join N$. The following hold:
\begin{enumerate}
\item $[0,N]_{\B}=0$;
\item $[M,N]_{\B}=[N,M]_{\B}$;
\item $[M,N]_{\B}\leq M\meet N$;
\item if $N\leq L$ then $[M,N]_{\B}\leq [M,L]_{\B}$ as subobjects of $A$;
\item $q_{J}[M,N]_{\B}\leq [q_{J}M,q_{J}N]_{\B}$;
\item $[M\times M',N\times N']_{\B}=[M,N]_{\B}\times [M',N']_{\B}$;
\item $q_{J}[M,N]_{\B}= [q_{J}M,q_{J}N]_{\B}$ as soon as $q_{J}(M\meet N)=q_{J}M\meet q_{J}N$, which happens, for instance, when either $M\leq N$ or $J\leq M\meet N$;
\item $[M,N]_{\B}$ is the smallest normal subobject $J$ of $M\join N$ such that $q_{J}M$ and $q_{J}N$ commute.
\end{enumerate}
\end{theorem}
\begin{proof}
The first property holds because, for any object $N$, the square
\[
\xymatrix@!0@=3.5em{N \ar[d]_-{q_{N}} \ar@{=}[r]^-{q_{0}} & N \ar[d] \\
0 \ar@{=}[r] & 0}
\]
is a double central extension with respect to $\B$. Property (2) follows from the symmetry of Diagram~\eqref{Square-Commutator}; see~\cite{EverHopf} for a detailed explanation. (3) follows from the definition of $[M,N]_{\B}$. To see this, consider the diagram
\[
\xymatrix@!0@C=7.5em@R=4em{K[[r_{0}]_{\B}] \ar@<.5ex>[d]^-{k_{1}} \ar@<-.5ex>[d]_-{k_{0}} \ar[r]^-{\ker [r_{0}]_{\B}} & [R_{M}\square R_{N}]_{\B} \ar@<.5ex>[r]^-{[r_{1}]_{\B}} \ar@<-.5ex>[r]_-{[r_{0}]_{\B}} \ar@<.5ex>[d]^-{[p_{1}]_{\B}} \ar@<-.5ex>[d]_-{[p_{0}]_{\B}} & [R_{N}]_{\B} \ar@<.5ex>[d]^-{[\pi_{1}]_{\B}} \ar@<-.5ex>[d]_-{[\pi_{0}]_{\B}}\\
K[[\rho_{0}]_{\B}] \ar[r]^-{\ker [\rho_{0}]_{\B}} \ar[d]_-{l} & [R_{M}]_{\B} \ar@<.5ex>[r]^-{[\rho_{1}]_{\B}} \ar@<-.5ex>[r]_-{[\rho_{0}]_{\B}} \ar[d]_-{\mu_{R_{M}}} & [M\join N]_{\B} \ar[d]^-{\mu_{M\join N}} \\
M \ar[r]_-{\ker \rho_{0}} & R_{M} \ar@<.5ex>[r]^-{\rho_{1}} \ar@<-.5ex>[r]_-{\rho_{0}} & M\join N.}
\]
Since $[M,N]_{\B}$, being the kernel of $\langle [p_{0}]_{\B},[r_{0}]_{\B}\rangle$, may be computed as the meet of the kernels of $[p_{0}]_{\B}$ and $[r_{0}]_{\B}$, it is also the kernel of~$k_{0}$. Hence, considered as a subobject of $M\join N$ via Remark~\ref{Remark-Subobject}, it is a subobject of $M$ through the morphism $l\comp k_{1}$. Likewise, $[M,N]_{\B}$ is contained in $N$.

The fourth property follows from the functoriality of the construction of $[-,-]_{\B}$. So does the fifth. To see that the relative commutator preserves binary products, it suffices to note that the zero-dimensional commutator~$[-]_{\B}$ preserves them, and that joins commute with products. The former property is well known. It is a consequence of the fact that the reflector $I\colon \A\to\B$ preserves pullbacks of split epimorphisms along split epimorphisms (because the components of the unit are extensions) together with the fact that a split epimorphism of split epimorphisms in $\Ext\A$ is always a three-fold extension. The latter property holds because the product of two regular pushouts is a regular pushout: products of pullbacks are pullbacks, products of regular epis are regular epis.

To prove (7), first of all recall that the square~\eqref{Diagram-unitsquare} induced by the unit~$\eta$ is a pushout of regular epimorphisms for any regular epimorphism $f$, by the Birkhoff condition. Hence, by Lemma~\ref{Rotlemma}, the zero-dimensional commutator $[-]_{\B}\colon \A\to\A$ preserves extensions. 

Now assume that $q_{J}(M\meet N)=q_{J}M\meet q_{J}N$. Then by Lemma~\ref{Lemma-Three-Fold-Extension} the left hand side commutative cube
\[
\xymatrix@!0@=3.5em{& \tfrac{M\join N}{J} \ar[rr] \ar@{.>}[dd] && \tfrac{M\join N}{M\join J} \ar[dd] \\
M\join N \ar[rr] \ar[dd] \ar[ru] && \tfrac{M\join N}{M} \ar[dd] \ar[ru]_-{\beta}\\
& \tfrac{M\join N}{N\join J} \ar@{.>}[rr] && 0\\
\tfrac{M\join N}{N} \ar[rr] \ar@{.>}[ru]^-{\alpha} && 0 \ar@{=}[ru] }
\qquad
\xymatrix@!0@=3.5em{& R_{\tfrac{M\join J}{J}}\square R_{\tfrac{N\join J}{J} } \ar@<.5ex>[rr] \ar@<-.5ex>[rr] \ar@{.>}@<.5ex>[dd] \ar@{.>}@<-.5ex>[dd] && R_{\tfrac{N\join J}{J}} \ar@<.5ex>[dd] \ar@<-.5ex>[dd]\\
R_{M}\square R_{N} \ar@<.5ex>[rr] \ar@<-.5ex>[rr] \ar@<.5ex>[dd] \ar@<-.5ex>[dd]\ar[ru] && R_{N} \ar@<.5ex>[dd] \ar@<-.5ex>[dd] \ar[ru]\\
& R_{\tfrac{M\join J}{J} } \ar@{.>}@<.5ex>[rr] \ar@{.>}@<-.5ex>[rr] && \tfrac{M\vee N}{J} \\
R_{M} \ar@<.5ex>[rr] \ar@<-.5ex>[rr] \ar@{.>}[ru] && M\join N \ar[ru] }
\]
is a three-fold extension. As a consequence, so are all the commutative cubes in the right hand side diagram, being pullbacks of three-fold extensions. This is still true if we apply the functor $[-]_{\B}$ to the right hand side diagram, since $[-]_{\B}$ preserves extensions and because a split epimorphism of extensions is a double extension, and a split epimorphism of double extensions a three-fold extension. The identity in (7) now follows. 

If~${M\leq N}$ then $q_{J}(M\meet N)=q_{J}M=q_{J}M\meet q_{J}N$. If, on the other hand, we assume that $J\leq M\meet N$, then the morphism $\alpha$ and, by symmetry, also $\beta$, are isomorphisms. This implies that the left hand side cube above is a three-fold extension, so that $q_{J}(M\meet N)=q_{J}M\meet q_{J}N$ by Lemma~\ref{Lemma-Three-Fold-Extension}.

Properties (3) and (7) together imply that $q_{[M,N]_{\B}}M$ and $q_{[M,N]_{\B}}N$ commute. Using (5) it is now easily seen that $[M,N]_{\B}$ is minimal amongst all $J$ such that ${[q_{J}M,q_{J}N]_{\B}=0}$.
\end{proof}

It was shown in~\cite{EverVdL4} that two normal subobjects of an $\Omega$-group commute in the sense of~\cite{EveraertCommutator} if and only if they commute in the sense of our Definition~\ref{Definition-Commuting-Subobjects}. Since both notions of relative commutator satisfy the same universal property (see Theorem~\ref{Proposition-Commutator-Properties} (8)), we find:

\begin{proposition}\label{Theorem-Characterisation}
Let $\A$ be a variety of $\Omega$-groups and $\B$ a subvariety of $\A$. Given any two normal subobjects $M$ and $N$ of an object $A$ of $\A$, we have
\[
[M,N]^{\Omega}_{\B}\cong [M,N]_{\B}.
\]
In particular, the commutator $[M,N]^{\Omega}_{\B}$ is zero if and only if the double extension~\eqref{Double-Extension-MN} is central.\noproof
\end{proposition}

\begin{remark}\label{Remark-examples}
This already gives us the examples worked out in~\cite{EveraertCommutator}: precrossed modules vs.\ crossed modules, where the relative commutator is the Peiffer commutator, for instance. An example which is not a consequence of this theorem---loops vs.\ groups, where the relative commutator is an associator---was considered in the article~\cite{EverVdL4}. Another example which falls outside the scope of~\cite{EveraertCommutator} is the case of compact Hausdorff topological groups vs.\ profinite groups. Here, the relative commutator $[M,N]_{\B}$ is the connected component of the intersection $M\cap N$, as follows from results in~\cite{Everaert-Gran-TT}. More generally, in any situation where the reflector $I\colon \A\to\B$ is \emph{protoadditive}~\cite{EG-honfg,Everaert-Gran-TT} (for instance, when $\A$ is abelian), one has the identity $[M,N]_{\B}=[M\cap N]_{\B}$ for any object $A$ of $\A$ and any pair of normal subobjects~$M$ and~$N$ of~$A$. 

The ``absolute'' case of abelianisation is treated in the following section.
\end{remark}

\begin{remark}\label{Remark-Preservation}
It suffices to consider the case $\B=0$ (where $0$ is the category with one object and one arrow) to see that the equality in Statement~(5) of Theorem~\ref{Proposition-Commutator-Properties} does not hold in general. The case $\B=0$ shows, furthermore, that unlike the Smith/Pedicchio commutator---cf.\ Lemma~\ref{Proposition-Smith-Commutator-Properties}---the commutator $[-,-]_{\B}$ need not preserve binary joins.
\end{remark}

\section{The absolute case: abelianisation}\label{Section-Huq}
In the case of $\Omega$-groups, the relative commutator $[-,-]_{\B}^{\Omega}$ in~$\A$ reduces to the Higgins commutator when $\B$ is the Birkhoff subcategory~$\Ab\A$ of all abelian objects of $\A$. Likewise, when $\A$ is an arbitrary semi-abelian category and $\B$ is~$\Ab\A$, the relative $[-,-]_{\B}$ is nothing but the Huq commutator. To show this we take a detour via the Smith/Pedicchio commutator of equivalence relations. First, in Lemma~\ref{Proposition-Commutator-Is-Smith}, we prove that the equivalence relation corresponding to the commutator of two normal subobjects is exactly the commutator of the equivalence relations corresponding to those normal subobjects. Then we prove Proposition~\ref{Proposition-Huq-Is-Smith} which states that the Huq commutator of a pair of normal subobjects $M$ and $N$ of an object~$A$ is the normalisation of the Smith/Pedicchio commutator of the corresponding equivalence relations, when $M\join N=A$. Combining both results, we obtain Theorem~\ref{Proposition-Commutator-Is-Huq}: given any two normal subobjects~$M$ and~$N$ of $A$, their Huq commutator $[M,N]^{\Huq}$, computed in $M\join N$, coincides with $[M,N]_{\Ab\A}$.

\subsection{The commutator of equivalence relations}
In his book~\cite{Smith}, Smith introduced a commutator of equivalence relations in the context of Mal'tsev varieties. It was extended to a purely categorical setting by Pedicchio~\cite{Pedicchio} and may be presented in a manner which is similar to the definition of the Huq commutator of normal subobjects~\cite{Borceux-Bourn, BG}.

Let $A$ be an object of a semi-abelian category $\A$. The largest equivalence relation on $A$ is denoted by $\nabla_{A}=(A\times A,\pi_{0},\pi_{1})$ and the smallest one by~$\Delta_{A}=(A,1_{A},1_{A})$.

Two equivalence relations $R=(R,r_{0},r_{1})$ and $S=(S,s_{0},s_{1})$ on~$A$ are said to \defn{centralise each other} when they admit a \defn{centralising double relation}
\begin{equation}\label{Centralising-Relation}
\vcenter{\xymatrix@!0@=3.5em{C \ar@<.5ex>[r] \ar@<-.5ex>[r] \ar@<.5ex>[d] \ar@<-.5ex>[d] & S \ar@<.5ex>[d] \ar@<-.5ex>[d]\\
R \ar@<.5ex>[r] \ar@<-.5ex>[r] & A,}}
\end{equation}
i.e., a (unique) double equivalence relation $C$ on $R$ and $S$ such that any of the four commutative squares in~\eqref{Centralising-Relation} is a pullback. (Then all of the commutative squares in~\eqref{Centralising-Relation} are pullbacks.) $R$ and $S$ centralise each other if and only if there exists a partial Mal'tsev operation on $R$ and $S$, a morphism $p\colon{R\times_{A}S\to A}$ which satisfies $p(\alpha,\alpha,\gamma)=\gamma$ and $p(\alpha,\gamma,\gamma)=\alpha$. 

The \defn{commutator} $[R,S]^{\Smith}$ of $R$ and $S$ is the universal equivalence relation on $A$ which, when divided out, makes them centralise each other. Consider the pullback
\[
\vcenter{\xymatrix@!0@=4em{R\times_{A}S \ar@<.5ex>[d]^-{p_{R}} \ar@<.5ex>[r]^-{p_{S}} \ar@{}[rd]|<<{\pullback} & S \ar@<.5ex>[d]^-{s_{0}} \ar@<.5ex>[l]^-{i_{S}}\\
R \ar@<.5ex>[r]^-{r_{1}} \ar@<.5ex>[u]^-{i_{R}} & A \ar@<.5ex>[l] \ar@<.5ex>[u]}}
\]
of $r_{1}$ and $s_{0}$; then $[R,S]^{\Smith}$ is the kernel pair $R[q]$ of the morphism $q$ in the diagram
\[
\xymatrix@!0@=3.5em{& R \ar[ld]_-{i_{R}} \ar@{.>}[d] \ar[rd]^-{r_{0}}\\
R\times_{A}S \ar@{.>}[r] & Q & A \ar@{.>}[l]|-{q}\\
& S \ar[lu]^-{i_{S}} \ar[ru]_-{s_{0}} \ar@{.>}[u]}
\]
where the dotted arrows denote the colimit of the outer square. The direct images~$qR$ and~$qS$ of~$R$ and~$S$ along the regular epimorphism $q$ centralise each other; hence $R$ and $S$ do so if and only if $[R,S]^{\Smith}=\Delta_{A}$.

The following properties of this commutator will be useful for us.

\begin{lemma}\label{Proposition-Smith-Commutator-Properties}\cite{Borceux-Bourn, Bourn-Gran-Maltsev, Pedicchio}
Let $R$, $S$, $S'$ be equivalence relations on an object~$A$ and $f\colon{A\to B}$ a regular epimorphism. The following hold:
\begin{enumerate}
\item $[\Delta_{A},S]^{\Smith}=\Delta_{A}$;
\item $[R,S]^{\Smith}=[S,R]^{\Smith}$;
\item $[R,S]^{\Smith}\leq R\meet S$;
\item if $S\leq S'$ then $[R,S]^{\Smith}\leq [R,S']^{\Smith}$;
\item $[R,S\join S']^{\Smith}=[R,S]^{\Smith}\join [R,S']^{\Smith}$;
\item if $[R,S]^{\Smith}=\Delta_{A}$ then $[fR,fS]^{\Smith}=\Delta_{B}$.\noproof
\end{enumerate}
\end{lemma}

The double central extensions with respect to the Birkhoff subcategory~$\Ab\A$ of abelian objects in a semi-abelian category $\A$ have been characterised in terms of this commutator of equivalence relations as follows. 

\begin{lemma}\cite{RVdL, EverVdL3}\label{Proposition-Central-With-Commutators} A double extension~\eqref{Diagram-Square} in a semi-abelian category~$\A$ satisfies
\[
[R[d],R[c]]^{\Smith}=\Delta_A=[R[d]\meet R[c], \nabla_A]^{\Smith}
\]
if and only if it is central with respect to $\Ab\A$.\noproof
\end{lemma}

This immediately implies that $[-,-]_{\Ab\A}$ corresponds to $[-,-]^{\Smith}$ in the following sense:

\begin{lemma}\label{Proposition-Commutator-Is-Smith}
Given any two normal subobjects $M$ and $N$ of $A$,
\[
[R_{M},R_{N}]^{\Smith}=R_{[M,N]_{\Ab\A}}.
\]
\end{lemma}
\begin{proof}
By definition, $M$ and $N$ commute when the square~\eqref{Double-Extension-MN} is a double central extension with respect to $\Ab\A$. According to Lemma~\ref{Proposition-Central-With-Commutators}, this happens if and only if
\begin{equation}\label{Two-Equalities}
[R_{M},R_{N}]^{\Smith}=\Delta_{M\join N}=[R_{M}\meet R_{N}, \nabla_{M\join N}]^{\Smith}.
\end{equation}
Using $\nabla_{M\join N}=R_{M}\join R_{N}$ we see that
\[
[R_{M}\meet R_{N}, \nabla_{M\join N}]^{\Smith}= [R_{M}\meet R_{N}, R_{M}]^{\Smith}\join [R_{M}\meet R_{N}, R_{N}]^{\Smith} \leq [R_{M},R_{N}]^{\Smith}
\]
and the second equality in~\eqref{Two-Equalities} follows from the first. Hence $[M,N]_{\Ab\A}$ is zero if and only if $[R_{M},R_{N}]^{\Smith}=\Delta_{M\join N}$. The commutator $[R_{M},R_{N}]^{\Smith}$ now coincides with $R_{[M,N]_{\Ab\A}}$ because these two equivalence relations satisfy the same universal property.
\end{proof}

\subsection{The Huq commutator} 
It is well known that in general, the Huq commutator does \emph{not} correspond to the commutator of equivalence relations: the relation~$R^{A}_{[M,N]^{\Huq}}$ need not be isomorphic to $[R_{M}^{A},R_{N}^{A}]^{\Smith}$ for arbitrary normal subobjects $M$ and $N$ of an object~$A$---a counterexample is given in~\cite{Bourn2004} for digroups, a variety of $\Omega$-groups. There are essentially two ways to remedy this situation. On the one hand, the context may be strengthened to that of \emph{Moore categories} by imposing the \emph{strong protomodularity} axiom~\cite{Borceux-Bourn, Rodelo:Moore}; but then the theory no longer applies to all varieties of $\Omega$-groups. On the other hand, it is known that the induced notions of centrality coincide in \emph{any} semi-abelian category (see~\cite[Proposition~2.2]{Gran-VdL}). That is to say, $R_{[M,N]^{\Huq}}^{A}$ is isomorphic to $[R_{M}^{A},R_{N}^{A}]^{\Smith}$ when $N$ is equal to $A$. In fact, according to an unpublished result by M.~Gran and the first author (presented here as Proposition~\ref{Proposition-Huq-Is-Smith} below) this assumption is too strong: as we shall see, the commutators coincide as soon as $A=M\join N$. 

Two coterminal morphisms $m\colon {M\to A}$ and $n\colon {N\to A}$ \defn{commute} when there exists a (necessarily unique) morphism $\varphi_{m,n}\colon{M\times N\to A}$ such that
\[
m=\varphi_{m,n}\comp \langle 1_{M},0\rangle\qquad\text{and}\qquad n=\varphi_{m,n}\comp \langle 0,1_{N}\rangle.
\]
It is clear that $m$ and $n$ commute if and only if their Huq commutator
\[
[m,n]^{\Huq}\colon [M,N]^{\Huq}\to A
\]
is zero, see Subsection~\ref{Subsection-Huq-Commutator}.

\begin{proposition}\label{Proposition-Huq-Is-Smith}
Given any two normal subobjects $M$ and $N$ of $A$ such that ${M\join N=A}$ we have $R_{[M,N]^{\Huq}}=[R_{M},R_{N}]^{\Smith}$.
\end{proposition}
\begin{proof}
We show that the representing normal monomorphisms $m$ and $n$ of~$M$ and~$N$ commute if and only if the equivalence relations $R_{M}$ and~$R_{N}$ centralise each other; the result then follows, because the commutators $[-,-]^{\Huq}$ and $[-,-]^{\Smith}$ satisfy the same universal property. One implication is Proposition~3.2 in~\cite{BG} which states that $m$ and $n$ commute whenever  $[R_{M},R_{N}]^{\Smith}$ is~$\Delta_{A}$. Indeed, if $p\colon{R_{M}\times_{A}R_{N}\to A}$ is a partial Mal'tsev operation on $R_{M}$ and $R_{N}$, then its restriction to $M\times N$ is the needed $\varphi_{m,n}$. 

To prove the other implication, suppose that $\varphi_{m,n}\colon{M\times N\to A}$ exists. By assumption, the morphism $\langle\begin{smallmatrix}m \\ n\end{smallmatrix}\rangle\colon {M+N\to A}$, and hence also $\varphi_{m,n}$, is a regular epimorphism. This implies that $R_{M}=\varphi_{m,n}(\varphi_{m,n}^{-1}R_{M})$ and $R_{N}=\varphi_{m,n}(\varphi_{m,n}^{-1}R_{N})$. Since the images of two equivalence relations which centralise each other still centralise each other (by (6) in Lemma~\ref{Proposition-Smith-Commutator-Properties}), it suffices to show that so do $\varphi_{m,n}^{-1}R_{M}$ and $\varphi_{m,n}^{-1}R_{N}$. Now these relations turn out to be particularly simple. Via Lemma~\ref{Lemma-Pullback}, the Noether isomorphism $N/(M\meet N)\cong(M\join N)/M$ implies that the left hand side square in the diagram with exact rows
\[
\xymatrix@!0@R=3.5em@C=6.5em{0 \ar[r] & M\times (M\meet N) \ar[d] \ar[r] & M\times N \ar[d]^-{\varphi_{m,n}} \ar[r] & \tfrac{N}{M\meet N} \ar[d]^-{\cong} \ar[r] & 0\\
0 \ar[r] & M \ar[r] & M\join N \ar[r] & \tfrac{M\join N}{M} \ar[r] & 0}
\]
is a pullback, so $\varphi_{m,n}^{-1}R_{M}=\nabla_{M}\times R_{M\meet N}^{N}$. Similarly, $\varphi_{m,n}^{-1}R_{N}=R_{M\meet N}^{M}\times \nabla_{N}$. Since
\[
[M\meet N,M]^{\Huq}=0=[M\meet N,N]^{\Huq},
\]
Proposition~2.2 in~\cite{Gran-VdL} may be used to see that both $[\nabla_{M},R_{M\meet N}^{M}]^{\Smith}=\Delta_{M}$ and $[R_{M\meet N}^{N},\nabla_N]^{\Smith}=\Delta_N$, so that $[\varphi_{m,n}^{-1}R_{M},\varphi_{m,n}^{-1}R_{N}]^{\Smith}=\Delta_{M\times N}$---which finishes the proof.
\end{proof}

Combining Lemma~\ref{Proposition-Commutator-Is-Smith} with Proposition~\ref{Proposition-Huq-Is-Smith}, we obtain

\begin{theorem}\label{Proposition-Commutator-Is-Huq}
Given any two normal subobjects $M$ and $N$ of $A$, their Huq commutator $[M,N]^{\Huq}$, computed in $M\join N$, coincides with $[M,N]_{\Ab\A}$.\noproof
\end{theorem}

\begin{remark}
Given any monomorphism $i\colon{A\to B}$, two coterminal morphisms $m\colon{M\to A}$ and $n\colon{N\to A}$ commute if and only if $i\comp m$ and $i\comp n$ commute---both in Huq's sense and relative with respect to any $\B$. This implies that the concept of ``commuting subobjects'' is independent of the surrounding object~$A$. As a consequence,
\[
\overline{[M,N]_{\Ab\A}}^{A}=[M,N]^{\Huq}.
\]
\end{remark}

\section{Further remarks}\label{Further-Remarks}

\subsection{Finding the right context}
We have defined the relative commutator in the framework of semi-abelian categories. However, looking at the diagram in the introduction, this is not entirely satisfactory, because: 
\begin{itemize}
\item Central extensions were defined in~\cite{Janelidze-Kelly} in the context of exact categories~$\A$, relative to a choice of admissible Birkhoff subcategory; and it was shown that if $\A$ is Mal'tsev (every reflexive relation internal in~$\A$ is an equivalence relation) then any Birkhoff subcategory is admissible. More recently, V.~Rossi proved in~\cite{Val} the admissibility of Birkhoff subcategories in a context which includes every regular Mal'tsev category that is ``almost exact'' in the sense that every regular epimorphism is an effective descent morphism.
\item The Huq commutator can be considered in a context,  as general as that of finitely cocomplete unital categories; in particular, in any finitely cocomplete pointed Mal'tsev category~\cite{Bourn-Huq}.
\end{itemize}
Thus one may ask if it is possible to consider the relative commutator in a more general context than that of semi-abelian categories, say, in finitely cocomplete, pointed, regular, ``almost exact'' Mal'tsev categories? We do not know the answer, but let us mention here two apparent obstacles and comment on either of these. 

(1) Double central extensions, on which concept the notion of relative commutator depends, were defined in~\cite{EGVdL} in the semi-abelian context. One reason for this was that the construction of the left adjoint to the inclusion functor $\CExt_{\B}\A\to \Ext\A$ given in~\cite{EGVdL} is only valid if $\A$ is semi-abelian (and~$\B$ is a Birkhoff subcategory of $\A$). In this case, the same construction can be applied to higher dimensions, giving us, in particular, a left adjoint to the inclusion functor $\CExt^2_{\B}\A\to \Ext^2\!\A$ of double central extensions into double extensions. The existence of the latter adjoint or, more precisely, of the reflection into $\CExt^2_{\B}\A$ of double extensions of the form \eqref{Double-Central-Extension} is what allows us to define the relative commutator. 

There is no a priori reason, though, why the left adjoints $\Ext\A\to\CExt_{\B}\A$ and $\Ext^2\!\A\to\CExt^2_{\B}\A$ could not exist when the category $\A$ is not semi-abelian. In fact, the former adjoint is known to exist in a wide variety of cases (see \cite{Janelidze-Kelly:Reflectiveness,Janelidze:Hopf}). For instance, it exists if $\A$ is a finitely cocomplete exact Mal'tsev category and $\B$ the Birkhoff subcategory of abelian objects, and in this case the characterisation of Lemma~\ref{Proposition-Central-With-Commutators} above remains valid (see~\cite{EverVdL3}).

(2) In an exact Mal'tsev category any pushout of regular epimorphisms is a regular pushout~\cite{Carboni-Kelly-Pedicchio}, and we have used this property to conclude the crucial fact that the square~\eqref{Double-Central-Extension} is always a double extension. Furthermore, we know from~\cite{Carboni-Kelly-Pedicchio} that in every regular, but not exact, Mal'tsev category there exist pushout squares of regular epimorphisms that are not double extensions. This seems to indicate that exactness is unavoidable in defining a relative commutator. However, we can say the following.
 
First of all we recall from~\cite{Bourn1996} that a finitely complete category $\A$ is Mal'tsev if and only if for any square of split epimorphisms
\[
\vcenter{\xymatrix@!0@=3.5em{X \ar@<.5ex>[d]^-{d} \ar@<-.5ex>[r]_-{c} & C \ar@<.5ex>[d]^-{g} \ar@<-.5ex>[l]\\
D \ar@<-.5ex>[r]_-{f} \ar@<.5ex>[u] & Z \ar@<-.5ex>[l] \ar@<.5ex>[u]}}
\]
which ``reasonably'' commutes (in the sense that it represents a split epimorphism in the category of split epimorphisms, with given splitting, in $\A$), the factorisation $\langle d,c\rangle\colon X\to D\times_Z C$ to the pullback of $f$ with $g$ is a strong epimorphism. A finitely complete pointed category $\A$ is called \defn{unital} if the same property holds, but only in the case where $Z$ is the zero object. Equivalently, $\A$ is unital if for any two objects $C$ and $D$ the ``product injections'' $\langle 0, 1_{C}\rangle\colon C\to D\times C$ and $\langle 1_{D},0\rangle\colon D\to D\times C$ are jointly strongly epimorphic~\cite{Bourn1996,Bourn2002}. A third characterisation of unital categories is given by the following proposition.

\begin{proposition}
If $\A$ is a finitely complete pointed category, then the first condition implies the second:
\begin{enumerate}
\item $\A$ is unital;
\item for any pair of strong epimorphisms $c$ and $d$
\[
\xymatrix@!0@=3.5em{
D & X \ar[l]_-d \ar[r]^-c & C}
\]
such that the kernels $\ker d$ and $\ker c$ are jointly strongly epimorphic, the induced morphism to the product $\langle d,c\rangle\colon X\to D\times C$ is a strong epimorphism.
\end{enumerate}
If, moreover, $\A$ has finite coproducts, then the two conditions are equivalent.
\end{proposition}
\begin{proof}
Assume that $\A$ is unital and that $d$ and $c$ are as in (2). First of all note that a morphism is a strong epimorphism if it is jointly strongly epimorphic with a zero morphism. Since $\ker d$ and $\ker c$ are jointly strongly epimorphic, and $d$ is a strong epimorphism, this implies that the composite $d\comp \ker c$ is strongly epimorphic. Similarly, $c\comp \ker d$ is a strong epimorphism. Since~$\A$ is unital, the product injections $\langle 0, 1_{C}\rangle$ and $\langle 1_{D},0\rangle$ are jointly strongly epimorphic, hence, by the above, so are $\langle 0, 1_{C}\rangle\comp c\comp \ker d=\langle d,c\rangle\comp \ker d$ and $\langle 1_{D},0\rangle\comp c\comp \ker c=\langle d,c\rangle\comp \ker c$. Hence $\langle d,c\rangle$ is a strong epimorphism. 

Conversely, for any two objects $D$ and $C$ of $\A$, applying condition (2) to the ``coproduct projections''
\[
\xymatrix@!0@=5em{
D & D+C \ar[l]_-{\langle\begin{smallmatrix}1_{D} \\ 0\end{smallmatrix}\rangle} \ar[r]^-{\langle\begin{smallmatrix}0 \\ 1_{C}\end{smallmatrix}\rangle} & C}
\]
gives us that the product injections $\langle 1_{D},0\rangle$ and $\langle 0, 1_{C}\rangle$ are jointly strongly epimorphic. Hence $\A$ is unital.
\end{proof}
Now suppose that $\A$ is finitely cocomplete, regular and unital. Then, in particular, any two normal subobjects $M$ and $N$ of an object $A$ in $\A$ admit a union~${M\cup N}$, and the above proposition implies that the square
\begin{equation}\label{Diagram-Non-Exact}
\vcenter{\xymatrix@!0@R=3.5em@C=5em{M\cup N \ar[r]^-{q_{M}} \ar[d]_-{q_{N}} & \tfrac{M\cup N}{M} \ar[d]\\
\tfrac{M\cup N}{N} \ar[r] & 0}}
\end{equation}
is a double extension (here $q_M$ and $q_N$ are the cokernels of the inclusions in~$M\cup N$ of $M$ and $N$, respectively). This indicates that it might be possible, after all, to consider the relative commutator in a non-exact context, but we would need to have an appropriate notion of double central extension. In that case, we could say that $M$ and $N$ commute if and only if the double extension~\eqref{Diagram-Non-Exact} is central.

\subsection{Stability under regular images}
We proved in Theorem~\ref{Proposition-Commutator-Properties} that 
\begin{equation}\label{Identity-Stable-Images}
 p[M,N]_{\B}=[pM,pN]_{\B}
\end{equation}
for any regular epimorphism $p\colon A\to B$ and normal subobjects $M$ and $N$ of~$A$ such that either $K[p]\leq M\meet N$ or $M\leq N$. As noted in Remark~\ref{Remark-Preservation}, this identity need not hold for arbitrary $p$, $M$ and $N$. However, we know from~\cite{Huq} that~\eqref{Identity-Stable-Images} \emph{does} hold for arbitrary $p$, $M$ and $N$ if $\B=\Ab\A$, and the same is true, for instance, for the Peiffer commutator of precrossed modules or the associator of loops (considered in~\cite{EverVdL4}). 

This suggests to look for necessary and sufficient conditions on the Birkhoff subcategory $\B$ for $[-,-]_{\B}$ to be \emph{stable under regular images}, i.e., for the identity~\eqref{Identity-Stable-Images} to hold for any regular epimorphism $p\colon A\to B$ and any normal subobjects $M$ and $N$ of $A$. We do not have a satisfactory answer to this question, although a characterisation of such $\B$ in the case of $\Omega$-groups was given in~\cite{EveraertCommutator}, in terms of the identities that define the subvariety $\B$. 

Let us just recall here the following necessary condition, again taken from the article~\cite{EveraertCommutator}: we need the subcategory $\Ab\A$ of abelian objects of~$\A$ to be contained in~$\B$. Indeed, if we assume that the relative commutator $[-,-]_{\B}$ is stable under regular images, and that $A$ is an abelian object with ``multiplication'' $\pi\colon {A\times A\to A}$, then
\begin{align*}
[A,A]_{\B} &= \Bigl[\pi \bigl(A\times 0\bigr), \pi \bigl(0 \times A\bigr)\Bigr]_{\B} = \pi \Bigl(\bigl[A \times 0, 0 \times A\bigr]_{\B}\Bigr)\\
 &\subseteq \pi\Bigl(\bigl(A \times 0\bigr) \meet \bigl(0 \times A\bigr)\Bigr) = 0.
\end{align*}
However, the converse is not true. The condition $\B\supseteq\Ab\A$ does \emph{not} imply the stability under regular images of $[-,-]_{\B}$; a counterexample was given in~\cite{EveraertCommutator}.

A similar question may be asked with respect to preservation of joins, see Remark~\ref{Remark-Preservation}.

\subsection{Higher dimensions}
In this article, we considered what we have called zero-dimensional, one-dimensional and two-dimensional relative commutators, but what about higher dimensions? Keeping in mind examples such as the associator of loops, this does not seem to be an unreasonable question to ask. Let us write $[L,M,N]_{\B}$ for a three-dimensional relative commutator defined on triples of normal subobjects $L$, $M$, $N$ of an object $A$ of $\A$, with respect to a Birkhoff subcategory $\B$ of $\A$. Then, for instance, if $\A$ is the variety of loops and $\B$ is the subvariety of groups, we would like that
\[
[L,M,N]_{\B}=[L,M,N],
\] 
where on the right hand side is the associator of loops considered in~\cite{EverVdL4}.

It is not clear to us what would be the appropriate definition of $n$-dimen\-sion\-al relative commutator (for $n\geq 3$), or whether it is even possible to obtain a convenient theory.

%

\end{document}